\documentclass{amsart}
\usepackage{amsfonts}
\usepackage{amsmath,amscd}
\usepackage{amssymb}
\usepackage{amsthm}
\usepackage{newlfont}

 \newtheorem{thm}{Theorem}[section]
 
 \newtheorem{lem}[thm]{Lemma}
 
 \theoremstyle{definition}
 
 \theoremstyle{remark}

 \numberwithin{equation}{section}

\begin{document}

\title[Generalized power graph of groups]
 {Generalized power graph of groups}

\author[A. Jafarzadeh, P. Niroomand  and M. Parvizi]{A. Jafarzadeh $^1$,  P. Niroomand$^2$  and M. Parvizi $^3$}

\address{$^{1}$ Department of Pure Mathematics, Ferdowsi University of Mashhad, Mashhad, Iran}
\email{jafarzadeh@um.ac.ir}

\address{$^{2}$ School of Mathematics and Computer Science,
Damghan University, Damghan, Iran}
\email{niromand@du.ac.ir, p$\_$niroomand@yahoo.com}

\address{$^{3}$ Department of Pure Mathematics, Ferdowsi University of Mashhad, Mashhad, Iran}
\email{parvizi@math.um.ac.ir}




\date{\today}


\begin{abstract}
The power graph of an arbitrary group $G$ is a simple graph with all elements of $G$ as its vertices and two vertices are adjacent if one is a positive power of another.  In this paper, we generalize this concept to a graph whose vertices are all elements of $G$ that generate a proper subgroup of $G$ and two elements are adjacent if the cyclic subgroups generated by which have non-trivial intersections. We concentrate on completeness and planarity of this graph.
\end{abstract}

\maketitle

\section{Introduction}

The investigation of graphs which come from algebraic structures is a growing area of research in mathematics.
firstly, Erd\"{o}s  introduced a graph
associated with a group, called its commuting graph, and asked whether its clique number is finite or not when all of its cliques are finite \cite{Erdos}. Later, other graphs were introduced based on algebraic structures, including the prime graph of a group \cite{Williams}, the zero-divisor graph of a ring \cite{Beck}, and the power graph of a semigroup or a group \cite{Cameron,ChakrabartyGhosh, Doostabadi}. 
The power graph of a group $G$ is a simple undirected graph $\mathcal{P}(G)$ that consists of all elements of $G$ as its vertices, and two vertices are adjacent if one is a positive power of another.

In this paper, we intend to generalize the notion of power graph of groups and then we obtain some information and facts related to this concept.   For a group $G$, let $\mathcal{GP}(G)$, the generalized power graph of $G$, be the graph which has all elements of $G$ that generate a proper subgroup of $G$ as its vertices and two elements are adjacent if the cyclic subgroups generated by which have non-trivial intersections. With a slight difference in vertices one can consider $\mathcal{P}(G)$ as a subgraph of $\mathcal{GP}(G)$ and the later is a subgraph of the intersection graph, of course by a different view of the vertices. The main subject of this paper, is to study this graph, mainly the completeness and planarity are discussed. In the case of abelian groups, results are more accurate and explicit. Despite abelian groups, non-abelian groups have their own difficulties to deal with. However, we try to obtain some valuable results on completeness and planarity of $\mathcal{GP}(G)$. In section 2 and 3, we focus on completeness of the generalized power graph of  abelian and non-abelian groups, respectively. Sections 4 and 5 are devoted to the planarity of this graph.

\section{Completeness of the generalized power graph of abelian groups}
To study the completeness of the mentioned graph, the orders of the elements of $G$ play an important role in choosing the techniques and methods. Let $G$ be an abelian group which is not necessarily finite. We wish to determine when $\mathcal{GP}(G)$ is complete. We have to consider torsion and torsion-free groups separately. Trivially, no mixed group, i.e. a group which is nether torsion, nor torsion-free, has complete graph. So the above two cases cover the entire subject. Starting with torsion groups we have the following very useful lemma:

\begin{lem}\label{complete-torsion}
Let $G$ be a torsion group and $\mathcal{GP}(G)$ be a complete graph. Then $G$ is a $p$-group for some prime $p$.
\end{lem}

\begin{proof}
The existence of elements of composite orders, implies the existence of elements of different prime power orders which contradicts to the completeness of $\mathcal{GP}(G)$.
\end{proof}

Now the following theorem for torsion abelian groups classifies all groups with the complete graph:

\begin{thm}
Let $G$ be a torsion abelian group, then $\mathcal{GP}(G)$ is complete if and only if $G\cong \mathbb{Z}_{p^n}$ or $G\cong \mathbb{Z}_{p^{\infty}}$, for some prime $p$. In fact, the first case occurs exactly when $G$ is bounded and the second case when $G$ is unbounded.
\end{thm}

\begin{proof}
Trivially, both groups $\mathbb{Z}_{p^n}$ and $\mathbb{Z}_{p^{\infty}}$ have complete graphs. Conversely, Lemma \ref{complete-torsion} shows that $G$ must be a $p$-group and completeness of $\mathcal{GP}(G)$ forces $G$ to have only one subgroup of order $p$, so the result holds.
\end{proof}

For torsion-free abelian groups we have the following simple structure:

\begin{thm}
Let $G$ be a torsion-free abelian group. Then $\mathcal{GP}(G)$ is complete if and only if $G$ is a subgroup of $\mathbb{Q}$, the additive group of rational numbers.
\end{thm}

\begin{proof}
Trivially $\mathbb{Q}$ and its subgroups has the mentioned property (note that these are the only locally cyclic abelian groups). Conversely, let $G$ be a torsion-free abelian group with complete graph, then for all $a,b\in G$, the subgroup generated by $a$ and $b$ must be an infinite cyclic group and hence $G$ is locally cyclic as desired.
\end{proof}

\section{Completeness of the generalized power graph of non-abelian groups}
In the case of non-abelian groups, there are not explicit results in infinite case; but in finite case, we have the foregoing nice result:

\begin{thm}\label{finite-non-abelian}
Let $G$ be a finite non-abelian group, then $\mathcal{GP}(G)$ is complete if and only if $G\cong Q_{2^n}$ for some $n\in \mathbb{N}$.
\end{thm}

\begin{proof}
The completeness of $\mathcal{GP}(Q_{2^n})$ is evidently seen because of the unique subgroup of order $p$. Conversely, Lemma \ref{complete-torsion} shows that $G$ is a $p$-group and by completeness it has only one subgroup of order $p$, so $G$ is a Quaternion group, as claimed.
\end{proof}

For infinite non-abelian groups we have the following partial result:

\begin{thm}
Let $G$ be a nilpotent torsion-free group with complete graph, then $G$ is abelian.
\end{thm}

\begin{proof}
Let $G$ be nilpotent of class $c>1$ and $a_1,\ldots, a_c\in G$ be arbitrary elements. If $[a_1,\ldots,a_{c-1}]\neq 1$ then $[a_1,\ldots,a_{c-1}]^r=a_1^s$ for some integers $r$ and $s$. This implies $1=[a_1,\ldots,a_{c-1},a_c^s]=[a_1,\ldots,a_{c-1},a_c]^s$. As $G$ is torsion-free, we have $\gamma_{c}(G)=1$, a contradiction. Hence the result holds.
\end{proof}

The following theorem classifies torsion non-abelian groups whose graphs are complete with the extra condition of being locally finite:

\begin{thm}
Let $G$ be a torsion non-abelian group with complete graph. If $G$ is locally finite, then $G$ is a (possibly infinite) generalized Quaternion group.
\end{thm}

\begin{proof}
Let $H$ be a finitely generated subgroup of $G$. By Theorem \ref{finite-non-abelian}, $H$ is a Quaternion group and hence  $G$ is a generalized Quaternion group.
\end{proof}

At the end of this section, the following result on the generalized power graph of $p$-groups shows not only the complete graphs belong to $p$-groups, but the generalized power graph of a $p$-group is not far from of being complete:

\begin{thm}
Let $G$ be a $p$-group. Then every component of $\mathcal{GP}(G)$ is complete. Furthermore, the number of its components is exactly the number of distinct subgroups of order $p$ in $G$.
\end{thm}

\begin{proof}
We simply prove $\mathrm{dim}~G=1$, the second part is trivial. Let $\mathrm{dim}~G=n>1$ and $a\leftrightarrow x_1\leftrightarrow \ldots \leftrightarrow x_{n-1}\leftrightarrow b$ be a path of length $n$. $\langle a\rangle \cap \langle x_1\rangle$ and $\langle x_1\rangle \cap \langle x_2\rangle$ are non-trivial subgroups of the cyclic group $\langle x_1\rangle$, and hence they have non-trivial intersection. So $a\leftrightarrow x_2$ which implies $n=1$.
\end{proof}

\section{Planarity of the generalized power graph of abelian groups}
In this section, we concentrate on abelian groups and study the planarity of $\mathcal{GP}(G)$ in this case. The techniques used here are only for finite abelian groups. We show that the groups with planar graphs have small exponents. In fact, ignoring elementary abelian $p$-groups, only $\mathbb{Z}_4$ and $\mathbb{Z}_6$ have planar graphs. We do this job in a sequence of lemmas.

\begin{lem}\label{4primes}
Let $G$ be a finite abelian group whose order is divisible by at least four distinct primes, then $\mathcal{GP}(G)$ is not planar.
\end{lem}

\begin{proof}
Let $p$, $q$, $r$ and $s$ be four primes dividing $|G|$. By choosing elements of orders $pq, pr, qr, pqs$ and $prs$ in a cyclic group of order $pqrs$, one can see that $\mathcal{GP}(G)$ has a subgraph isomorphic to $K_5$ and hence it is not planar.
\end{proof}

Lemma \ref{4primes} shows that the number of prime divisors of the order of an abelian group $G$ with planar generalized power graph does not exceeds 3. The following lemma shows the prime divisors are not too large, too.

\begin{lem}\label{p>7}
Let $G$ be a finite  group with planar generalized power graph. Then dividing the order of $G$ by $p$ implies  $p\leq 5$.
\end{lem}

\begin{proof}
It is easy to see that the graph of $\mathbb{Z}_p$ for $p\geq 7$ is not planar, so if $p$ divides the order of $G$, then we must have $p\leq 5$.
\end{proof}

The above lemmas show that for finding finite abelian groups with planar generalized power graph we have to consider only $\pi$-groups with $\pi \subseteq \{2,3,5\}$. The next lemma refines this search in some sense.

\begin{lem}
Let $G$ be a finite abelian group with planar graph. Then either $G$ is a $p$-group for some $p\in \{2,3,5\}$ or $G$ is a $\{2,3\}$-group.
\end{lem}

\begin{proof}
By the above-mentioned sentence, it is enough to show that each $\{2,5\}$ or $\{3,5\}$-group fails to have planar generalized power graph. Let $G$ be a $\{2,5\}$-group and $x$ be an element of order $10$ in $G$. Then $\{x,x^2,x^4,x^6,x^8\}$ induces the complete graph $K_5$ and hence $\mathcal{GP}(G)$ is not planar. The other case can be proved similarly.
\end{proof}

Now the following theorem completes the classification of finite abelian groups with planar generalized power graph.

\begin{thm}
Let $G$ be a finite abelian group. Then $\mathcal{GP}(G)$ is planar if and only if one of the following holds:

\begin{enumerate}
\item $G$ is an elementary abelian $2$-group.
\item $G$ is an elementary abelian $3$-group.
\item $G$ is an elementary abelian $5$-group.
\item $G\cong \mathbb{Z}_4$ 
\item $G\cong \mathbb{Z}_6$
\end{enumerate}
\end{thm}

\begin{proof}
It is easy to see that all groups listed in the theorem have planar graphs. Let $G$ be a finite abelian group with planar generalized power graph. As  mentioned before, $G$ must be a $p$-group for some $p\in \{2,3,5\}$ or it is a $\{2,3\}$-group. If $G$ is a $p$-group for $p\in \{3,5\}$ then its exponent must be $p$, otherwise $\mathcal{GP}(\mathbb{Z}_{p^2})$ is a non-planar subgraph of $\mathcal{GP}(G)$. If $G$ is a $2$-group, then with a similar argument we conclude that the exponent of $G$ divides 4. Now, if  the exponent of $G$ equals 4 and $G$ is not cyclic of order $4$, then it has a subgroup isomorphic to $\mathbb{Z}_4\oplus \mathbb{Z}_2$ which has non-planar graph. Finally, let $G$ be a $\{2,3\}$-group, with a similar argument, the exponent of $G$ is $6$ and because of the non-planarity of the graph of $\mathbb{Z}_6\oplus \mathbb{Z}_2$, $G$ must be cyclic of order $6$.
\end{proof}

\section{Planarity of the generalized power graph of non-abelian groups}
In this section, we consider non-abelian groups and try to classify all groups having planar generalized power graph. The results are partial and only covers $p$-groups, but some similarities may be found here. For instance, the prime divisors of the order of $G$ are  the same as the case of abelian groups.

We only consider the planarity of non-abelian $p$-groups for $p\in \{2,3,5\}$, the last two cases can be handled easily as follows.

\begin{thm}
Let $G$ be a non-abelian $p$-group with planar generalized power graph. If $p\in \{3,5\}$, then the exponent of $G$ equals $p$. Furthermore, if the order of $G$ is $p^n$, then $\mathcal{GP}(G)$ is a union of $\frac{p^n-1}{p-1}$ components each of which is isomorphic to $K_{p-1}$.  
\end{thm}

\begin{proof}
If the exponent of $G$ is not equal to $p$, then $G$ has a subgroup isomorphic to $\mathbb{Z}_{p^2}$ whose graph is not planar.  For the second part it is easy to see that the number of components is equal to the number of subgroups of order $p$ and each component is complete.
\end{proof}
  
As usual $p=2$ differs in conclusions. See the following theorem.

\begin{thm}
Let $G$ be a finite non-abelian $2$-group with planar graph then $G\cong D_8$. 
\end{thm}

\begin{proof}
Sinse $G$ is not abelian, the exponent of $G$ is greater than 2. If it is greater than 4, then $G$ has a subgroup of order $8$ which has a non-planar graph so the exponent of $G$ is 4. Let $x$ be an element of order $4$. Since $\mathbb{Z}_4\oplus \mathbb{Z}_2$ has non-planar graph we have $C_{G}(x)=\langle x \rangle$ and so $Z(G)\subseteq \langle x \rangle$ is of order $2$ or $4$. But $G$ is non-abelian so $Z(G)$ must be of order $2$. This implies the exponent of $\frac{G}{Z(G)})$ equals 2 and so $G$ is an extra special $2$-group. If $|G|>8$, then $Q_8\subseteq G$ so $|G|=8$ and hence $G\cong D_8$.
\end{proof}


\end{document}